\documentclass{amsart}
\newtheorem{Thm}{Theorem}[section]
\newtheorem{Prop}[Thm]{Proposition}
\newtheorem{Cor}[Thm]{Corollary}
\newtheorem{Lem}[Thm]{Lemma}
\newtheorem{Rem}[Thm]{Remark}
\numberwithin{equation}{section}
\begin{document}

\title[Graphs and the $p$-harmonic boundary]
{Graphs of bounded degree and the $p$-harmonic boundary}

\author[M. J. Puls]{Michael J. Puls}
\address{Department of Mathematics \\
John-Jay College-CUNY \\
445 West 59th Street \\
New York, NY 10019 \\
USA}
\email{mpuls@jjay.cuny.edu}

\begin{abstract}
Let $p$ be a real number greater than one and let $G$ be a connected graph of bounded degree. In this paper we introduce the $p$-harmonic boundary of $G$. We use this boundary to characterize the graphs $G$ for which the constant functions are the only $p$-harmonic functions on $G$. It is shown that any continuous function on the $p$-harmonic boundary of $G$ can be extended to a function that is $p$-harmonic on $G$. Some properties of this boundary that are preserved under rough-isometries are also given. Now let $\Gamma$ be a finitely generated group. As an application of our results we characterize the vanishing of the first reduced $\ell^p$-cohomology of $\Gamma$ in terms of the cardinality of its $p$-harmonic boundary. We also study the relationship between translation invariant linear functionals on a certain difference space of functions on $\Gamma$, the $p$-harmonic boundary of $\Gamma$ and the first reduced $\ell^p$-cohomology of $\Gamma$.
\end{abstract}


\date{September 8, 2010}
\maketitle

\section{Introduction}\label{Introduction}
Let $p$ be a real number greater than one and let $\Gamma$ be a finitely generated infinite group. There has been some work done relating various boundaries of $\Gamma$ and the nonvanishing of the first reduced $\ell^p$-cohomology space $\bar{H}_{(p)}^1(\Gamma)$ of $\Gamma$ (to be defined in Section \ref{lpcohomology}). It was shown in Chapter 8, section C2 of \cite{Gromov} (also see \cite{Elek} ) that if the $\ell^p$-corona of $\Gamma$ contains more than one element, then $\bar{H}^1_{(p)} (\Gamma) \neq 0$. In \cite{PulsADM} it was shown that if there is a Floyd boundary of $\Gamma$ containing more than two elements, and if the Floyd admissible function satisfies a certain decay condition, then $\bar{H}^1_{(p)} (\Gamma) \neq 0$. However, it is unknown if the converse of either of these two results is true. The motivation for this paper is to find a boundary for $\Gamma$ that characterizes the vanishing of $\bar{H}_{(p)}^1 (\Gamma)$ in terms of the cardinality of this boundary. In this paper we will show that the $p$-harmonic boundary (which we define in Section \ref{pharmbd}) does the trick. The reason this boundary gives the desired result is that $\bar{H}^1_{(p)}(\Gamma) = 0$ if and only if the only $p$-harmonic functions on $\Gamma$ are the constants, \cite[Theorem 3.5]{Puls}. We will show in Section \ref{lpcohomology} that the cardinality of the $p$-harmonic boundary is 0 or 1 if and only if the only $p$-harmonic functions on $\Gamma$ are the constants. Hence, $\bar{H}^1_{(p)}(\Gamma) = 0$ if and only if the cardinality of the $p$-harmonic boundary is $0$ or $1$.

$L_p$-cohomology was investigated first in \cite{Goldshtein} for the case of Riemannian manifolds. Gromov has studied $\ell^p$-cohomology for finitely generated groups, and in the more general setting of graphs with bounded degree, see Chapter 8 of \cite{Gromov}. In particular, it was shown in \cite{CheegerGromov} that the first reduced $\ell^2$-cohomology of a finitely generated amenable group is zero. It was conjectured by Gromov in \cite[page 150]{Gromov} that this result is also true for all real numbers $p>1$. This is our main justification for choosing to study the $p$-harmonic boundary in the discrete setting. If enough insight can be gained into this boundary, then we may be able to develop the tools needed to compute the $p$-harmonic boundary of a finitely generated amenable group. This of course would resolve Gromov's conjecture.

More information about the first reduced $L_p$-cohomology ($L_2$-cohomology) can be found in the papers \cite{PansuTorino, PansuPot, PansuHel, Tessera} for various manifolds, and in the papers \cite{BekkaValette, Bourdon, BMV, Elekcourse, MartinValette, PulsCMB, Puls, PulsADM} for finitely generated groups. As was implied earlier there is a strong connection between the vanishing of the first reduced $L_p$-cohomology ($L_2$-cohomology) and the nonexistence of nonconstant $p$-harmonic functions (2-harmonic functions), see \cite[Proposition 4.11]{Tessera} for a proof in the case of homogeneous Riemannian manifolds. Thus papers that study $p$-harmonic functions are useful in trying to determine if the first reduced $L_p$-cohomology vanishes. The papers \cite{Coulhon, Grigoryan} study $p$-harmonic functions on manifolds, and the papers \cite{H2, positivepharm, finiteenergy, Soardirough, Yamasaki1977} examine $p$-harmonic functions on graphs.

I would like to thank Peter Linnell for many useful comments on a preliminary version of this paper. I would also like to thank the referee for some suggestions that greatly improved the exposition of the paper. The research for this paper was partially supported by PSC-CUNY grant 60123-38 39.

\section{Definitions and statement of main results}\label{Outlinestatementmain}

In this section we define many of the concepts studied in this paper. We also give an outline of the paper and state our results. Let $p$ be a real number greater than one and let $\Gamma$ be a finitely generated infinite group.  It turns out that the definition of the $p$-harmonic boundary for $\Gamma$ does not depend on the group law of $\Gamma$. Due to this fact we are able to define this boundary in the more general setting of a graph. The reason why we can do this is that we can associate a graph, called the Cayley graph of $\Gamma$, with $\Gamma$. The vertex set for this graph are the elements of $\Gamma$, and $x_1, x_2 \in \Gamma$ are joined by an edge if and only if $x_1 = x_2s^{\pm 1}$ for a generator $s$ of $\Gamma$. We now proceed to define the desired boundary.

\subsection{The $p$-harmonic boundary}\label{pharmbd}

Let $G$ be a graph with vertex set $V_G$ and edge set $E_G$. We will write $V$ for $V_G$ and $E$ for $E_G$ if it is clear what the graph $G$ is. For $x \in V, \mbox{ deg}(x)$ will denote the number of neighbors of $x$ and $N_x$ will be the set of neighbors of $x$. A graph $G$ is said to be of {\em bounded degree} if there exists a positive integer $k$ such that $\mbox{deg}(x) \leq k$ for every $x \in V$. A path in $G$ is a sequence of vertices $x_1, x_2, \dots, x_n$ where $x_{i+1} \in N_{x_i}$ for $1 \leq i \leq n-1$. A graph $G$ is connected if any two given vertices of $G$ are joined by a path. All graphs considered in this paper will be countably infinite, connected, of bounded degree with no self-loops. Assign length one to each edge in $E_G$, then the graph $G$ is a metric space with respect to the shortest path metric. Let $d_G( \cdot, \cdot)$ denote this metric. So if $x, y \in V$, then $d_G(x, y)$ is the length of the shortest path joining $x$ and $y$. We will drop the subscript $G$ from $d_G( \cdot, \cdot)$ when it is clear what graph $G$ we are working with. Finally, if $x \in V$, then $B_n(x)$ will denote the metric ball that contains all elements of $V$ that have distance less than $n$ from $x$.

Let $G$ be a graph with vertex set $V$ and let $p$ be a real number greater than one. In order to construct the $p$-harmonic boundary of $G$ we need to first define the space of bounded $p$-Dirichlet finite functions on $G$. For any $S \subset V$, the outer boundary $\partial S$ of $S$ is the set of vertices in $V\setminus S$ with at least one neighbor in $S$. For a real-valued function $f$ on $S \cup \partial S$ we define the $p$-th power of the {\em gradient}, the {\em $p$-Dirichlet sum}, and the {\em $p$-Laplacian} of $x \in S$ by
\begin{equation*}
\begin{split}
 \vert Df (x) \vert^p  & = \sum_{y \in N_x} \vert f(y) - f(x) \vert^p,   \\
 I_p (f, S)            & = \sum_{x \in S} \vert Df (x) \vert^p, \\
  \Delta_p f (x)    & = \sum_{y \in N_x} \vert f(y) - f(x) \vert^{p-2} (f(y) - f(x)).
\end{split}
\end{equation*} 
In the case $1 < p < 2$, we make the convention that $\vert f(y) - f(x) \vert^{p-2} (f(y) - f(x)) = 0$ if $f(y) = f(x)$. Let $S \subseteq V$. A function $f$ is said to be $p$-harmonic on $S$ if $\Delta_p f(x) = 0$ for all $x \in S$. We shall say that $f$ is {\em $p$-Dirichlet finite} if $I_p(f,V) < \infty$. The set of all $p$-Dirichlet finite functions on $G$ will be denoted by $D_p(G)$. Under the following norm $D_p(G)$ is a reflexive Banach space,
$$ \parallel f \parallel_{D_p} = \left( I_p(f,V) + \vert f(o) \vert^p \right)^{1/p},$$
where $o$ is a fixed vertex of $G$ and $f \in D_p(G)$. Denote by $HD_p(G)$ the set of $p$-harmonic functions on $V$ that are contained in $D_p(G)$. Let $\ell^{\infty}(G)$ denote the set of bounded functions on $V$ and let $\parallel f \parallel_{\infty} = \sup_V \vert f \vert$ for $f \in \ell^{\infty}(G)$. Set $BD_p(G) = D_p(G) \cap \ell^{\infty}(G)$. The set $BD_p(G)$ is a Banach space under the norm 
$$\parallel f \parallel_{BD_p} = \left( I_p(f, V)\right)^{1/p} + \parallel f \parallel_{\infty},$$
where $f \in BD_p(G)$. Set $BHD_p(G) = HD_p(G) \cap BD_p(G).$ It turns out that $BD_p(G)$ is closed under pointwise multiplication. To see this let $f,h \in BD_p(G)$ and set $a = \sup_V |f|$ and $b = \sup_V |h|$. It follows from Minkowski's inequality that
\begin{equation}
 \left( I_p( fh, V)\right)^{1/p} \leq b\left( I_p(f, V)\right)^{1/p} + a \left( I_p (h,V)\right)^{1/p}. \label{eq:productbound}
\end{equation}
Thus $fh \in BD_p(G)$. Using the above inequality we obtain
\begin{equation*}
\begin{split}
 \parallel fh \parallel_{BD_p} & \leq \left( \left( I_p(f, V)\right)^{1/p} + a\right)\left(\left(I_p(h,V) \right)^{1/p} + b \right)\\
                               & = \parallel f \parallel_{BD_p} \parallel h \parallel_{BD_p}.
\end{split}
\end{equation*}
Hence $BD_p(G)$ is an abelian Banach algebra. A character on $BD_p(G)$ is a nonzero homomorphism from $BD_p(G)$ into the complex numbers. We denote by $Sp(BD_p(G))$ the set of characters on $BD_p(G)$. With respect to the weak $\ast$-topology, $Sp(BD_p(G))$ is a compact Hausdorff space. The space $Sp(BD_p(G))$ is known as the spectrum of $BD_p(G)$. Let $C(Sp(BD_p(G)))$ denote the set of continuous functions on $Sp(BD_p(G))$. For each $f \in BD_p(G)$ a continuous function $\hat{f}$ can be defined on $Sp(BD_p(G))$ by $\hat{f}(\tau) = \tau(f)$. The map $f \rightarrow \hat{f}$ is known as the Gelfand transform.

Define a map $i \colon V \rightarrow Sp(BD_p(G))$ by $(i(x))(f) = f(x)$. For $x \in V$, define $\delta_x$ by $\delta_x(v) = 0$ if $v \neq x$ and $\delta_x (x) = 1$. Let $x, y \in V$ and suppose $i(x) = i(y)$, then $(i(x))(\delta_x) = (i(y))(\delta_x)$ which implies $\delta_x(x) = \delta_x(y)$. Thus $i$ is an injection. If $f$ is a nonzero function in $BD_p(G)$, then there exists an $x \in V$ such that $\hat{f}(i(x)) \neq 0$ since $\hat{f}(i(x)) = f(x)$. Hence $BD_p(G)$ is semisimple. Theorem 4.6 on page 408 of \cite{TaylorLay} now tells us that $BD_p(G)$ is isomorphic to a subalgebra of $C(Sp(BD_p(G)))$ via the Gelfand transform. Since the Gelfand transform separates points of $Sp(BD_p(G))$ and the constant functions are contained in $BD_p(G)$, the Stone-Weierstrass Theorem yields that $BD_p(G)$ is dense in $C(Sp(BD_p(G)))$ with respect to the sup-norm. The following proposition shows that $i(V)$ is dense in $Sp(BD_p(G))$, see \cite[Proposition 1.1(ii)]{Elek}for the proof.

\begin{Prop} \label{dense}
The image of $V$ under $i$ is dense in $Sp(BD_p(G))$.
\end{Prop}

When the context is clear we will abuse notation and write $V$ for $i(V)$ and $x$ for $i(x)$, where $x \in V$. The compact Hausdorff space $Sp(BD_p(G)) \setminus V$ is known as the $p$-Royden boundary of $G$, which we will denote by $R_p(G)$. When $p =2$ this is simply known as the Royden boundary of $G$. Let $\mathbb{R}G$ be the set of  real-valued functions on $V$ with finite support and let $B(\overline{\mathbb{R}G})_{D_p} = (\overline{\mathbb{R}G})_{D_p} \cap \ell^{\infty}(G)$. Suppose $(f_n)$ is a sequence in $B(\overline{\mathbb{R}G})_{D_p}$ that converges to a bounded function $f$ in the $BD_p(G)$-norm. It follows from $\parallel f-f_n \parallel_{D_p} \leq \parallel f - f_n \parallel_{BD_p}$ that $f \in (\overline{\mathbb{R}G})_{D_p}$. Thus $B(\overline{\mathbb{R}G})_{D_p}$ is closed in $BD_p(G)$ with respect to the $BD_p(G)$-norm. We are now ready to define the main object of study for this paper.

The {\em $p$-harmonic boundary} of $G$ is the following subset of the $p$-Royden boundary
$$\partial_p(G) \colon= \{ x \in R_p(G) \mid \hat{f}(x) = 0 \mbox{ for all } f \in B(\overline{\mathbb{R}G})_{D_p}\}.$$
 When $p = 2$ the $p$-harmonic boundary is known as the harmonic boundary. Our definition of the $p$-harmonic boundary is a direct generalization of the definition for the harmonic boundary. A good reference concerning the Royden and harmonic boundaries of graphs is Chapter VI of \cite{Soardi}.

An important fact about $B(\overline{\mathbb{R}G})_{D_p}$ is that it is an ideal in $BD_p(G)$. To see this let $f \in B(\overline{\mathbb{R}G})_{D_p}$ and $h \in BD_p(G)$. We need to show that $fh \in B(\overline{\mathbb{R}G})_{D_p}$. We claim that there exists a sequence $(f_n)$ in $\mathbb{R}G$ that converges pointwise to $f$, for which there exists a constant $M$ with $\vert f_n(x) \vert \leq M$ for all $n$ and for all $x \in V$, and has $I_p(f_n, V)$ bounded. To see the claim let $(u_n)$ be a sequence in $\mathbb{R}G$ that converges to $f$ in $D_p(G)$ and let $M= \sup_{x \in V} \vert f(x) \vert$. Set $f_n = \max(\min(u_n,M),-M).$ The sequence $(f_n)$ satisfies the claim above since $I_p(u_n, V)$ is bounded and $I_p(f_n, V) \leq I_p(u_n, V)$. Moreover, $(f_nh)$ is a sequence in $\mathbb{R}G$ that converges pointwise to $fh$. By inequality (\ref{eq:productbound}) we see that 
\[ I_p(f_nh, V) \leq (b(I_p(f_n, V))^{1/p} + M(I_p(h,V))^{1/p})^p, \]
where $b = \sup_{x \in V} \vert h(x) \vert.$ Since $I_p(f_nh, V)$ is bounded, Theorem 10.6 from page 177 of \cite{TaylorLay} says, by passing to a subsequence if necessary, that $(f_nh)$ converges weakly to a function $\overline{fh}$. Due to $B(\overline{\mathbb{R}G})_{D_p}$ being closed it follows that $\overline{fh} \in B(\overline{\mathbb{R}G})_{D_p}$. Because point evaluations by elements of $V$ are continuous linear functionals on $BD_p(G), (f_nh)$ also converges pointwise to $\overline{fh}$. Hence, $\overline{fh} = fh$ and $fh \in B(\overline{\mathbb{R}G})_{D_p}$.

\subsection{Statement of main results}\label{stmtresults}

In this subsection we outline the paper and state our main results. Recall that $p$ is a real number greater than one and that $o$ is a fixed vertex of $V$. We will use the notation $\#(A)$ to mean the cardinality of a set $A$ and $1_V$ will denote the function on $V$ that always takes the value one. Furthermore, $\ell^p(G)$ will be the set that consists of the functions on $V$ for which $\sum_{x \in V} \vert f(x) \vert^p < \infty$. The $\ell^p$-norm for $f \in \ell^p(G)$ is given by $\parallel f \parallel_p = (\sum_{x \in V} \vert f(x) \vert^p )^{\frac{1}{p}}.$ In Section \ref{pharmonicfunction} we give a quick review of some results about $p$-harmonic functions on graphs that will be needed in the sequel. In Section \ref{propertiespharmbound} we prove several results concerning $BD_p(G)$ and $\partial_p(G)$, including a characterization of when $\partial_p(G) = \emptyset$, when $BHD_p(G)$ consists precisely of the constant functions and a neighborhood base is given for the topology on $\partial_p(G)$.

Before we state some of our main results we need to state a theorem that will allow us to classify graphs in a nice way. We start by giving the following definition. The {\em $p$-capacity} of a finite subset $A$ of $V$ is defined by
$$ \mbox{Cap}_p(A, \infty, V) = \inf_u I_p(u, V)$$
where the infimum is taken over all finitely supported functions $u$ on $V$ such that $u = 1$ on $A$. The following theorem, which is Theorem 3.1 of \cite{Yamasaki1977}, will allow us to classify a graph $G$ in terms of the $p$-capacity of a finite set. See the reference for a proof of the theorem.
\begin{Thm} \label{yam77}
Let $A$ be a finite, nonempty subset of $V$. Then $\mbox{Cap}_p(A, \infty, V) =0$ if and only if $1_V \in B(\overline{\mathbb{R}G})_{D_p}.$
\end{Thm}

The following corollary follows directly from the theorem.
\begin{Cor}
Let $A$ and $B$ be nonempty finite subsets of $V$. Then $\mbox{Cap}_p (A, \infty, V) = 0$ if and only if $\mbox{Cap}_p (B, \infty, V) = 0$.
\end{Cor}
We shall say that a graph $G$ is {\em $p$-parabolic} if there exists a finite subset $A$ of $V$ such that $\mbox{Cap}_p(A, \infty, V) = 0$. If $G$ is not $p$-parabolic, then we shall say that $G$ is {\em $p$-hyperbolic}. Note that if $G$ is $p$-hyperbolic, then $\mbox{Cap}_p(A, \infty, V) > 0$ for all finite subsets $A$ of $V$. 

In Section \ref{Proofsmain} we will prove the following result, which is \cite[Theorem 4.6]{Soardi} for the case $p=2$. This result also generalizes Theorem 4.2 of \cite{positivepharm}.
\begin{Thm} \label{pparaconstant} Let $p$ be a real number greater than one and let $G$ be a graph. If $G$ is $p$-parabolic, then the constant functions are the only $p$-harmonic functions on $G$.
\end{Thm}
Identify the constant functions on $V$ by $\mathbb{R}$. By combining this theorem with Lemma 4.4 of \cite{H2} and Theorem \ref{plious} we obtain the following Liouville type theorem for $p$-harmonic functions
\begin{Thm} \label{constantchar}
Let $p$ be a real number greater than one. Then $HD_p(G) = \mathbb{R}$ if and only if the cardinality of $\partial_p(G)$ is either zero or one.
\end{Thm}
We also prove in Section \ref{Proofsmain}
\begin{Thm} \label{extendcont}
Let $p$ be a real number greater than one and let $G$ be a graph. If $f$ is a continuous function on $\partial_p(G)$, then there exists a $p$-harmonic function $h$ on $V$ such that $\lim_{n \rightarrow \infty} h(x_n) = f(x)$, where $x \in \partial_p(G)$ and $(x_n)$ is any sequence in $V$ that converges to $x$.
\end{Thm}
By combining the above theorem with the maximum principle and Corollary \ref{boundval} we obtain the following corollary, which is a generalization of both Theorem 4.3 of \cite{positivepharm} and Theorem 1.1 of \cite{finiteenergy}.
\begin{Cor} \label{finitedim}
Let $p$ be a real number greater than one and let $G$ be a graph. Assume that the $p$-harmonic boundary of $G$ contains a finite number of points, say $\{ x_1, x_2, \dots, x_n\}$. Then given real numbers $a_1, a_2, \dots, a_n \in \mathbb{R}$, there exists a bounded $p$-harmonic function $h$ that satisfies
\begin{equation} 
h(x_i) = a_i \mbox{ for } i = 1, 2, \dots , n.   \label{eq:existpharm}
\end{equation}
Conversely, each bounded $p$-harmonic function is uniquely determined by its values in (\ref{eq:existpharm}).
\end{Cor}

Let $(X, d_X)$ and $(Y, d_Y)$ be metric spaces. A map $\phi \colon X \rightarrow Y$ is said to be a {\em rough isometry} if it satisfies the following two conditions:
\begin{enumerate}
\item There exists constants $ a \geq 1, b \geq 0$ such that for $x_1, x_2 \in X$
$$ \frac{1}{a} d_X (x_1, x_2) - b \leq d_Y (\phi(x_1), \phi(x_2)) \leq a d_X(x_1, x_2) + b.  $$
\item There exists a positive constant $c$ such that for each $y \in Y$, there exists an $x \in X$ that satisfies $d_Y ( \phi(x), y) < c$.
\end{enumerate}

For a rough isometry $\phi$ there exists a rough isometry $\psi \colon Y \rightarrow X$ such that if $x \in X, y \in Y$, then $d_X( (\psi \circ \phi)(x), x) \leq a(c+b)$ and $d_Y((\phi \circ \psi)(y), y) \leq c$. The map $\psi$, which is not unique, is said to be a rough inverse for $\phi$. Whenever we refer to a rough inverse to a rough isometry in this paper, it will always satisfy the above conditions. In section \ref{roughproofs} we will prove
\begin{Thm} \label{homeomorphic}
Let $p$ be a real number greater than one and let $G$ and $H$ be graphs. If there is a rough isometry from $G$ to $H$, then $\partial_p(G)$ is homeomorphic to $\partial_p(H)$.
\end{Thm}
We will finish Section \ref{roughproofs} by proving
\begin{Thm} \label{bijection}
Let $p$ be a real number greater than one and let $G$ and $H$ be graphs. If there is a rough isometry from $G$ to $H$, then there is a bijection from $BHD_p(G)$ to $BHD_p(H)$.
\end{Thm}

The main result of \cite{Soardirough} is that if $G$ and $H$ are roughly isometric graphs then $HD_p(G) = \mathbb{R}$ if and only if $HD_p(H) = \mathbb{R}$. By Lemma 4.4 of \cite{H2} this is equivalent to $BHD_p(G) = \mathbb{R}$ if and only if $BHD_p(H) = \mathbb{R}$. Both Theorem \ref{homeomorphic} and Theorem \ref{bijection} are generalizations of this result.

We now return to the case of a finitely generated group $\Gamma$. In Section \ref{lpcohomology} we define the first reduced $\ell^p$-cohomology space $\bar{H}^1_{(p)} (\Gamma)$ of $\Gamma$. Then we will use our results on $p$-harmonic boundaries to prove
\begin{Thm} \label{lpchar}
Let $1 < p \in \mathbb{R}$. Then $\bar{H}^1_{(p)} (\Gamma) \neq 0$ if and only if $\#(\partial_p(\Gamma)) >1$.
\end{Thm}
It appears that there are not many explicit examples of the $p$-Royden boundary $R_p(G)$ for a given graph $G$. The only example that we know of is the paper \cite{Wysoczanski}, where the author gave an explicit description of $R_2(\mathbb{Z})$. We will conclude Section \ref{lpcohomology} by using Theorem \ref{lpchar} to compute the $p$-harmonic boundary for the case $\Gamma = \mathbb{Z}^n$. We will also compute the $p$-Royden boundary of nonamenable groups with infinite center, and for $\Gamma_1 \times \Gamma_2 \times \cdots \times \Gamma_n$, where $n \geq 2$, each $F_i$ is finitely generated and at least one of the $\Gamma_i$ is nonamenable. 

Let $E$ be a normed space of functions on a finitely generated group $\Gamma$. Let $f \in E$ and let $x \in \Gamma$. The right translation of $f$ by $x$, denoted by $f_x$, is the function $f_x(g) = f(gx^{-1})$, where $g \in \Gamma.$ Assume that if $f \in E$ then $f_x \in E$ for all $x\in \Gamma$; that is, that $E$ is right translation invariant. For the rest of this paper translation invariant will mean right translation invariant. We shall say that $T$ is a translation invariant linear functional (TILF) on $E$ if $T(f_x) = T(f)$ for $f \in E$ and $x \in \Gamma$. We will use TILFs to denote translation invariant linear functionals. A common question to ask is that if $T$ is a TILF on $E$, then is $T$ continuous? For background about the problem of automatic continuity see \cite{ Meisters, Saeki, Willis2, Woodward}. Define
$$ \mbox{Diff}(E) := \mbox{ linear span}\{f_x - f \mid f \in E, x \in \Gamma \}. $$
It is clear that $\mbox{Diff}(E)$ is contained in the kernel of any TILF on $E$. In Section \ref{tilf} we study TILFs on $D_p(\Gamma)/\mathbb{R}$. In particular we prove

\begin{Thm} \label{chartilf} 
Let $\Gamma$ be a finitely generated infinite group and let $1 < p \in \mathbb{R}$. Then $\#(\partial_p(\Gamma)) > 1$ if and only if there exists a nonzero continuous TILF on $D_p(\Gamma)/\mathbb{R}.$
\end{Thm}

It was shown in \cite{Willis1} that if $\Gamma$ is nonamenable, then the only TILF on $\ell^p(\Gamma)$ is the zero functional. (Consequently every TILF is automatically continuous!). We will conclude Section \ref{tilf} by showing that this result is not true for $D_p(\Gamma)/\mathbb{R}$.

\section{Review of $p$-harmonic functions}\label{pharmonicfunction}
In this section we will give some facts about $p$-harmonic functions on graphs that will be needed in this paper. Most of this material is from Section 3 of \cite{H2}, where a more comprehensive treatment, including proofs, is given.\\

{\bf Existence:} Let $S$ be a finite subset of $V$. For any function $f$ on $\partial S$, there exists an unique function $h$ on $S \cup \partial S$ which is $p$-harmonic on $S$ and equals $f$ on $\partial S$. \\
Moreover, in the proof of existence, it was shown that the $p$-harmonic function $h$ satisfies $\mbox{min}_{y \in \partial S} f(y) \leq h(x) \leq \mbox{max}_{y \in \partial S} f(y)$ for all $x \in S$. \\

{\bf Minimizer property:} Let $h$ be a $p$-harmonic function on a finite subset $S$ of $V$. Then $I_p(h, S) \leq I_p(f, S)$ for all functions $f$ on $S \cup \partial S$ that satisfy $f = h$ on $\partial S$. \\

{\bf Convergence:} Let $(S_n)$ be an increasing sequence of finite connected subsets of $V$ and let $U = \cup_i S_i$. Let $(h_i)$ be a sequence of functions on $U \cup \partial U$, such that $h_i(x) \rightarrow h(x) < \infty$ for every $x \in U \cup \partial U$. If $h_i$ is $p$-harmonic on $S_i$ for all $i$, then $h$ is $p$-harmonic on $U$. \\

{\bf Comparison principle:} Let $h$ and $u$ be $p$-harmonic functions on a finite subset $S$ of $V$. If $h \geq u$ on $\partial S$, then $ h \geq u$ on $S$.
\\

We conclude this section by proving the maximum principle for bounded $p$-harmonic functions on $V$.
\begin{Lem} \label{maxonV}
Let $h$ be a $p$-harmonic function on $V$. If there exists an $x \in V$ such that $h(x) \geq h(y)$ for all $y \in V$, then $h$ is constant on $V$.
\end{Lem} 
\begin{proof}
Let $x \in V$ such that $h(x) \geq h(x')$ for all $x' \in V$. Because $\sum_{y \in N_x} \vert h(y) - h(x) \vert^{p-2} h(y) = \sum_{y \in N_x} \vert h(y) - h(x) \vert^{p-2} h(x)$ we see that $h(x) = h(y)$ for all $y \in N_x$. Thus $h(x) = h(z)$ for all $ z \in V$ since $G$ is connected. 
\end{proof}

\section{Preliminary results}\label{propertiespharmbound}
In this section we will give some results about $\partial_p(G)$ and $BD_p(G)$ that will be needed in the sequel. Most of the results given from Proposition \ref{noelement} through through \ref{vanishbound} are given in the first two sections of Chapter VI in \cite{Soardi} for the case of $p=2$. However, our presentation and some of our proofs are different. Recall that $o$ is a fixed vertex of the graph $G$. We begin with
\begin{Lem}\label{infinitelength}
If $x \in \partial_p(G)$ and $(x_n)$ is a sequence in $V$ that converges to $x$, then $d(o, x_n) \rightarrow \infty$ as $n \rightarrow \infty.$
\end{Lem}
\begin{proof}
Let $x \in \partial_p(G)$ and suppose $(x_n) \rightarrow x$, where $(x_n)$ is a sequence in $V$. Let $B$ be a positive real number. Define a function $\chi_B$ on $V$ by $\chi_B (y) = 1$ if $d(o, y) \leq B$ and $\chi_B(y) = 0$ if $d(o, y) > B$. Since $\chi_B$ has finite support it is an element of $\mathbb{R}G$. Suppose there exists a real number $M$ such that $d(o, x_n) \leq M$ for all $n$. Then $\widehat{\chi_M} (x) = \lim_{n \rightarrow \infty} \chi_M (x_n) = 1$, a contradiction. Thus $d(o, x_n) \rightarrow \infty$ as $n \rightarrow \infty.$
\end{proof}
We now characterize $p$-parabolic graphs in terms of $\partial_p(G)$.
\begin{Prop} \label{noelement}
Let $G$ be a graph and let $1 < p \in \mathbb{R}$. Then $\partial_p(G) = \emptyset$ if and only if $G$ is $p$-parabolic.
\end{Prop}
\begin{proof}
Assume $G$ is $p$-parabolic and suppose $\partial_p(G)\neq \emptyset$. Let $x \in \partial_p(G)$ and let $(x_n)$ be a sequence in $V$ that converges to $x$. Then $\widehat{1_V} (x) = \lim_{n \rightarrow \infty} \widehat{1_V} (x_n) =1$. By Theorem \ref{yam77}, $1_V \in B(\overline{\mathbb{R}G})_{D_p}$, which says that $\widehat{1_V}(x) = 0$, a contradiction. Hence if $G$ is $p$-parabolic, then $\partial_p(G) = \emptyset$. 

Now suppose that $G$ is $p$-hyperbolic. Then $1_V \notin B(\overline{\mathbb{R}G})_{D_p}.$ Since $B(\overline{\mathbb{R}G})_{D_p}$ is an ideal in the commutative ring $BD_p(G)$ there exists a maximal ideal $M$ in $BD_p(G)$ containing $B(\overline{\mathbb{R}G})_{D_p}$. Using the correspondence between maximal ideals in $BD_p(G)$ and $ Sp(BD_p(G))$ there is an $x \in Sp(BD_p(G))$ that satisfies $Ker(x) = M$. So $\hat{f}(x) = x(f)=0$ for all $f \in B(\overline{\mathbb{R}G})_{D_p}$. For each $y \in V$ there exists an $f \in \mathbb{R}G$ (in particular $\delta_y$) such that $y(f) = f(y) \neq 0$ which means that $x$ cannot be in $V$. Furthermore, if $x \in R_p(G) \setminus \partial_p(G)$, then there exists an $f \in B(\overline{\mathbb{R}G})_{D_p}$ for which $\hat{f}(x) \neq 0$. However, this implies that $B(\overline{\mathbb{R}G})_{D_p}$ is not contained in $M$. Therefore $x \in \partial_p(G)$.
\end{proof}
For the rest of this paper, unless otherwise stated, we will assume that $1_V \notin B(\overline{\mathbb{R}G})_{D_p}$; that is, $G$ is $p$-hyperbolic.

Let $f$ and $h$ be elements in $BD_p(G)$ and let $1 < p \in \mathbb{R}$. Define 
$$ \langle \triangle_p h, f \rangle \colon = \sum_{x \in V} \sum_{y \in N_x} \vert h(y) - h(x) \vert^{p-2} (h(y) - h(x))(f(y) - f(x)).$$
This sum exists since $\sum_{x \in V} \sum_{y \in N_x} \vert \vert h(y) - h(x) \vert^{p-2} (h(y) - h(x))\vert^q = I_p(h, V) < \infty$ where $\frac{1}{p} + \frac{1}{q} = 1$. The next few lemmas will be used to help show the uniqueness of the decomposition of $BD_p(G)$ that will be given in Theorem \ref{decomp}.
\begin{Lem} \label{diffconst}
Let $f_1$ and $f_2$ be functions in $D_p(G)$. Then $\langle \triangle_p f_1 - \triangle_p f_2, f_1 - f_2 \rangle = 0$ if and only if $f_1 - f_2$ is constant on $V$.
\end{Lem}
\begin{proof}
Let $f_1, f_2 \in D_p(G)$ and assume there exists an $x \in V$ with a $y \in N_x$ such that $f_1(x) - f_1(y) \neq f_2(x) - f_2(y)$. Define a function $f \colon [0,1] \rightarrow \mathbb{R}$ by
$$ f(t) = \sum_{x \in V} \sum_{y \in N_x} \vert f_1(y) - f_1(x) + t (( f_2(y) - f_2(x)) - (f_1(y) -f_1(x))) \vert^p.$$
Observe that $f(0) = I(f_1, V)$ and $f(1) = I(f_2, V)$. A derivative calculation gives
$$ f'(0) = p \langle \triangle_p f_1, f_2 - f_1 \rangle = -p \langle \triangle_p f_1, f_1 - f_2 \rangle.$$
It follows from Proposition 5.4 of \cite{Ekeland} that $I_p (f_2, V) > I_p( f_1, V) - p\langle \triangle_p f_1, f_1 - f_2 \rangle$. Similarly, $I_p(f_1, V) > I_p(f_2, V) - p \langle \triangle_p f_2, f_2 - f_1 \rangle$. Hence, $p\langle \triangle_p f_1 - \triangle_p f_2, f_1-f_2 \rangle > 0$ if there exists an $x \in V$ with $y \in N_x$ that satisfies $f_1(x) - f_1(y) \neq f_2(x) - f_2(y)$.

Conversely, suppose $f_1 - f_2$ is constant on $V$. We immediately see that $\langle \triangle_p f_1 - \triangle_pf_2, f_1 - f_2 \rangle =0.$
\end{proof}
\begin{Lem} \label{diracpharm}
Let $h \in BD_p(G)$. Then $h \in BHD_p(G)$ if and only if $\langle \triangle_p h, \delta_x \rangle = 0$ for all $x \in V$.
\end{Lem}
\begin{proof}
Let $x \in V$ and let $h \in BD_p(G)$. The lemma follows from 
$$\langle \triangle_p h, \delta_x \rangle = - 2 (\mbox{deg}(x))\sum_{y\in N_x} \vert h(x) - h(y)\vert^{p-2} (h(y) - h(x)).$$
\end{proof}
\begin{Rem}
The lemma implies that if $h \in BHD_p(G)$, then $\langle \triangle_p h, f \rangle = 0$ for all $f \in \mathbb{R}G.$
\end{Rem}
\begin{Lem} \label{killsubspace}
If $h \in BHD_p(G)$ and $f \in B(\overline{\ell^p(G)})_{D_p}$, then $\langle \triangle_p h, f \rangle = 0$.
\end{Lem}
\begin{proof}
Let $f \in B(\overline{\ell^p (G)})_{D_p}$ and let $h \in BHD_p(G)$. Then there exists a sequence $(f_n)$ in $\mathbb{R}G$ such that $\parallel f - f_n \parallel_{D_p} \rightarrow 0$ as $n \rightarrow \infty$ since $(\overline{\mathbb{R}G})_{D_p} = (\overline{\ell^p(G)})_{D_p}$. Now
\begin{eqnarray*} 
0  & \leq  &  \vert \langle \triangle_p h, f \rangle \vert = \vert \langle \triangle_p h, f-f_n \rangle \vert  \\
   & =   & \left| \sum_{x \in V} \sum_{y \in N_x} \vert h(y) - h(x) \vert^{p-2} (h(y) - h(x) ) ((f - f_n)(x) - (f-f_n)(y))\right|  \\ 
   & \leq & \sum_{x \in V} \sum_{y \in N_x} \vert h(y) - h(x) \vert^{p-1} \vert (f - f_n)(x) - (f - f_n)(y) \vert  \\
   & \leq & \left( \sum_{x \in V} \sum_{y \in N_x} (\vert h(y) - h(x) \vert^{p-1} )^q \right)^{1/q} \left( I_p (f - f_n, V) \right)^{1/p} \rightarrow 0
\end{eqnarray*}
as $n \rightarrow \infty$. The last inequality follows from H\"{o}lder's inequality.
\end{proof}
We now state Clarkson's Inequality, which will be needed in the proof of our next result.

Let $f_1$ and $f_2$ be elements of $D_p(G)$, if $2 \leq p \in \mathbb{R}$ then
\[ I_p(f_1+f_2)+I_p(f_1-f_2) \leq 2^{p-1}\left( I_p(f_1) + I_p(f_2) \right) \]
and for $1 < p \leq 2$
\[ \left( I_p(f_1 + f_2) \right)^{1/(p-1)} + \left( I_p(f_1 - f_2)\right)^{1/(p-1)} \leq 2 \left( I_p(f_1) + I_p(f_2)\right)^{1/(p-1)}. \]

We will now give a decomposition of $BD_p(G)$ that will be crucial in later work.
\begin{Thm} \label{decomp} 
Let $1 < p \in \mathbb{R}$ and suppose $f \in BD_p(G)$. Then there exists a unique $u \in B(\overline{\ell^p(G)})_{D_p}$ and a unique $h \in BHD_p(G)$ such that $f = u + h$.
\end{Thm}
\begin{proof} 
Remember our standing assumption that $1_V \notin B(\overline{\ell^p(G)})_{D_p}$. Let $f \in BD_p(G)$. Since $f$ is bounded there exists real numbers $a$ and $b$ for which $a \leq f(x) \leq b$ is satisfied by all $x \in V$. Denote by $h_n$ the function that is $p$-harmonic on $B_n(o)$ and equal to $f$ on $V\setminus B_n(o)$. Because $\min_{y \in \partial B_n(o)} f(y) \leq h_n(x) \leq \max_{y \in \partial B_n(o)} f(y)$ for all $x \in B_n(o)$, we have that $a \leq h_n \leq b$ for each $n \in \mathbb{N}$. Furthermore, if $m > n$, then $I_p(h_m) \leq I_p(h_n)$. Set $r_n = I_p(h_n)$ and denote the limit of the bounded decreasing sequence $(r_n)$ by $r$. We are still assuming that $m > n$. By the minimizing property of $p$-harmonic functions $I_p(h_m, V) \leq I_p( (h_n + h_m)/2, V )$ since $(h_n + h_m)/2 = h_m$ on $V\setminus B_m(o)$. Using Clarkson's Inequality we obtain the following: \\
For $ 2 \leq p \in \mathbb{R}$, \\
\begin{eqnarray*}
  r_m  &  \leq   &  I_p\left( \frac{h_n + h_m}{2}, V\right)  \\
       &  \leq   & I_p\left( \frac{h_n + h_m}{2}, V\right) + I_p \left(\frac{h_n - h_m}{2}, V\right) \\
       &  \leq   & 2^{p-1}\left( I_p\left(\frac{h_n}{2}, V\right) + I_p\left(\frac{h_m}{2}, V\right) \right) \\
       &   =     & \frac{1}{2} \left( I_p(h_n, V) + I_p(h_m, V) \right)
\end{eqnarray*}
and for $1 < p \leq 2$, \\
   \begin{eqnarray*}
     r_m^{1/(p-1)} & \leq  & \left( I_p\left( \frac{h_n+h_m}{2}, V \right) \right)^{1/(p-1)} \\
                   &  \leq  & \left( I_p\left( \frac{h_n+h_m}{2}, V \right) \right)^{1/(p-1)} + \left( I_p \left( \frac{h_n - h_m}{2},V\right)\right)^{1/(p-1)} \\
                   & \leq &  2 \left( I_p \left( \frac{h_n}{2}, V \right) + I_p \left( \frac{h_m}{2}, V \right)\right)^{1/(p-1)}.
   \end{eqnarray*}
Letting $m,n \rightarrow \infty$ it follows that $I_p(\frac{h_n+h_m}{2}, V) \rightarrow r$ and $I_p(\frac{h_n-h_m}{2}, V) \rightarrow 0$. Also, $( |h_n(o)|)$ is a bounded sequence, thus $(h_n)$ is a Cauchy sequence in $D_p(G)$. Set $h$ equal to the limit function of the sequence $(h_n)$ in $D_p(G)$. Because $(h_n)$ also converges pointwise to $h$ the convergence property says that $h$ is $p$-harmonic. Clearly, $a \leq h \leq b$ on $V$, so $h \in BHD_p(G)$. Let $u$ be the limit function in $D_p(G)$ of the Cauchy sequence $(f- h_n)$. Since $f - h_n \in \mathbb{R}G$ for each $n$, we see that $u \in B(\overline{\mathbb{R}G})_{D_p}$. Thus $f = u+h$.

We will now show that this decomposition is unique. Suppose $f = u_1 + h_1 = u_2 + h_2$, where $u_1, u_2 \in B(\overline{\ell^p(G)})_{D_p}$ and $h_1, h_2 \in BHD_p(G)$. Lemma \ref{killsubspace} says that $\langle \triangle_p h_1 - \triangle_p h_2, h_1 - h_2 \rangle = \langle \triangle_p h_1 - \triangle_p h_2, u_2 - u_1 \rangle = 0$ since $u_1 - u_2 \in B(\overline{\ell^p(G)})_{D_p}$. However, $u_1 - u_2 = 0$ since $1_V \notin B(\overline {\ell^p(G)})_{D_p}$.
\end{proof}

We are now ready to prove the following
\begin{Thm} \label{maxprinciple} (Maximum principle)
Let $h$ be a nonconstant function in $BHD_p(G)$ and suppose $a$ and $b$ are real numbers for which $a \leq \hat{h} \leq b$ on $\partial_p(G)$. Then $a < h < b$ on $V$.
\end{Thm}
\begin{proof}
Since $\hat{h}$ is continuous on the compact space $Sp(BD_p(G))$ there is a number $c >0$ such that $b - \hat{h} \geq -c$ on $Sp(BD_p(G))$. Let $\epsilon > 0$ and set $F_{\epsilon} = \{ x \in Sp(BD_p(G)) \mid b-h + \epsilon \leq 0 \}.$ In order to prove the theorem we will first show that there exists an $f \in B(\overline{\mathbb{R}G})_{D_p}$ with $\hat{f} = 1$ on $F_{\epsilon}$ and $0 \leq \hat{f} \leq 1$ on $Sp(BD_p(G)).$ This $f$ will yield the following inequality 
\begin{equation} 
  cf + b - h + \epsilon \geq 0 \mbox{ on } Sp(BD_p(G)). \label{eq:greaterthanzero}
\end{equation}
We will then show that $b-h+\epsilon \geq 0$ on $V$. Combining this with Lemma \ref{maxonV} and the assumption $h$ is nonconstant will give $h < b$ on $V$.

Observe that $F_{\epsilon} \cap \partial_p(G) = \emptyset$ and $F_{\epsilon}$ is a closed subset of $Sp(BD_p(G))$. For each $x\in F_{\epsilon}$ there exists an $f_x \in B( \overline{\mathbb{R}G})_{D_p}$ for which $\hat{f}_x(x) \neq 0$. Since $B(\overline{\mathbb{R}G})_{D_p}$ is an ideal we may assume that $f_x \geq 0$ on $V$ and $\hat{f}_x(x) > 0$. Let $U_x$ be a neighborhood of $x$ in $Sp(BD_p(G))$ that satisfies $f_x(y) > 0$ for all $y \in U_x$. By compactness there exists $x_1, \dots, x_n$ for which $F_{\epsilon} \subseteq \cup_{j=1}^n U_{x_j}$. Set $g = \sum_{j=1}^n f_{x_j}$ and let $\alpha = \inf\{ g(x) \mid x \in F_{\epsilon} \}$. Clearly $\alpha > 0$ and $g \in B( \overline{\mathbb{R}G})_{D_p}$. Now define a function $f$ on $Sp(BD_p(G))$ by $f = \min(1, \alpha^{-1}g)$. Note that $0 \leq \hat{f} \leq 1$ on $Sp(BD_p(G))$ and $\hat{f} = 1$ and $F_{\epsilon}$. We still need to show that $f \in B(\overline{\mathbb{R}G})_{D_p}$. Let $(g_n)$ be a sequence in $\mathbb{R}G$ that converges to $g$ in $D_p(G)$, so $I_p((g-g_n), V) \rightarrow 0$ has $n \rightarrow \infty$. Set $f_n = \min(1, \alpha^{-1}g_n)$. The sequence $(f_n)$ converges pointwise to $f$. Furthermore, by passing to a subsequence if necessary, $(f_n)$ converges weakly to a function $\bar{f}$ in $D_p(G)$ since $I_p  (f_n, V)$ is bounded. Clearly $\bar{f}$ is bounded, so $\bar{f} \in B( \overline{\mathbb{R}G})_{D_p}$. It is also true $(f_n)$ converges pointwise to $\bar{f}$ because point evaluations by elements of $V$ are continuous linear functionals on $BD_p(G)$. Hence, $\bar{f} = f$ and $f \in BD_p(G)$. Inequality (\ref{eq:greaterthanzero}) is now established.

Next we will show that $b-h+\epsilon \geq 0$ on $V$. Put $v_{\epsilon} = cf + b - h +\epsilon$ and denote by $h_n$ the unique function that is $p$-harmonic on $B_n(o)$ and agrees with $v_{\epsilon}$ on $V \setminus B_n(o)$. We claim that $h_n \geq 0$ on $B_n(o)$. Suppose otherwise, then there exists an $x \in B_n(o)$ for which $h_n(x) < 0$. Define a function $h_n^{\ast}$ by  \\
\[ h_n^{\ast} = \left\{ \begin{array}{cl}  
     v_{\epsilon}   &    \mbox{for } x \in V\setminus B_n(o) \\
    \max(h_n, 0)    &     \mbox{for } x \in B_n(o)  \end{array} \right. .\]
Now $I_p(h_n^{\ast}, B_n(o)) < I_p(h_n, B_n(o))$, but this contradicts the minimizer property for $p$-harmonic functions. This proves the claim. By using the argument used in the proof of Theorem \ref{decomp} we have that $(h_n)$ converges to a bounded $p$-harmonic function $\bar{h}$ and that there exists a $v \in B(\overline{\mathbb{R}G})_{D_p}$ such that $v_{\epsilon} = v + \bar{h}$. Furthermore $\bar{h} \geq 0$ on $V$ due to $h_n \geq 0$ for each $n$. The uniqueness part of Theorem \ref{decomp} says that $v = cf$ and $\bar{h} = b-h+\epsilon$. Hence $b \geq h - \epsilon$ on $V$. Thus $h < b$ on $V$.

Using a similar argument it can be shown that $a < h$ on $V$. Therefore, $a < h < b$ on $V$.
\end{proof}
We now characterize the functions in $BD_p(G)$ that vanish on $\partial_p(G)$. 
\begin{Thm} \label{vanishbound} 
Let $f \in BD_p(G)$. Then $f \in B(\overline{\ell^p(G)})_{D_p}$ if and only if $\hat{f}(x) = 0$ for all $ x \in \partial_p(G).$
\end{Thm}
\begin{proof}
Since $B(\overline{\ell^p(G)})_{D_p} = B(\overline{\mathbb{R}G})_{D_p}$ it follows immediately that $\hat{f}(x) =0$ for all $f \in B(\overline{\ell^p(G)})_{D_p}$ and all $x \in \partial_p(G)$.

Conversely, suppose $f \in BD_p(G)$ and $\hat{f} (x) = 0$ for all $x \in \partial_p(G)$. By Theorem \ref{decomp} we can write $f = u + h$, where $u \in B(\overline{\ell^p(G)})_{D_p}$ and $h \in BHD_p(G)$. Now $\hat{h} (x) = 0$ for all $x \in \partial_p(G)$ since $\hat{u}(x) = 0$. Therefore, $h=0$ by the maximum principle.
\end{proof}
As a consequence of the theorem we obtain
\begin{Cor} \label{boundval}
A function in $BHD_p(G)$ is uniquely determined by its values on $\partial_p(G)$.
\end{Cor}
\begin{proof}
Let $h_1$ and $h_2$ be elements of $BHD_p(G)$ with $\widehat{h_1} (x) = \widehat{h_2} (x)$ for all $x \in \partial_p(G)$. Then $h_1 - h_2 \in B(\overline{\ell^p(G)})_{D_p}$. Let $(f_n)$ be a sequence in $\ell^p(G)$ that converges to $h_1 - h_2$. Using Lemma \ref{killsubspace} we obtain $\langle \triangle_p h_1 - \triangle_p h_2, h_1 - h_2 \rangle = \lim_{n \rightarrow \infty} \langle \triangle_p h_1 - \triangle_p h_2, f_n \rangle = 0$. It now follows from Lemma \ref{diffconst} that $h_1 - h_2 =0$.
\end{proof}
We are now ready to give a characterization of when $BHD_p(G)$ is precisely the constant functions.
\begin{Thm} \label{plious}
Let $1 < p \in \mathbb{R}$. Then $BHD_p(G) \neq \mathbb{R}$ if and only if $\#(\partial_p(G)) > 1$.
\end{Thm}
\begin{proof}
Suppose that $\#(\partial_p (G)) = 1$ and that $x \in \partial_p(G)$. Let $h \in BHD_p(G)$. Then $\hat{h}(x) = c$ for some constant $c$. It follows from Corollary \ref{boundval} that the function $h(x) = c$ for all $x \in V$ is the only function in $BHD_p(G)$ with $\hat{h} (x) = c$. Hence $BHD_p(G) = \mathbb{R}$.

Conversely, suppose $\#(\partial_p(G)) > 1$. Let $x, y \in \partial_p(G)$ such that $x \neq y$ and pick an $f \in BD_p(G)$ that satisfies $x(f) \neq y(f)$. By Theorem \ref{vanishbound} $f \notin B(\overline{\ell^p(G)})_{D_p}$. It now follows from Theorem \ref{decomp} and Theorem \ref{vanishbound} that there exists an $h \in BHD_p(G)$ with $\hat{h}(z) = \hat{f}(z)$ for all $z \in \partial_p(G)$. Since $V$ is dense in $Sp(BD_p(G))$ there exists sequences $(x_n)$ and $(y_n)$ in $V$ such that $(x_n)(h) \rightarrow x(h)$ and $(y_n)(h) \rightarrow y(h)$. Hence $\lim_{n \rightarrow \infty} h(x_n) = x(h) \neq y(h) = \lim_{n \rightarrow \infty} h(y_n)$. Therefore, $h$ is not constant on $V$.
\end{proof}
Before we give our next results we need to define the important concept of a $D_p$-massive subset of a graph. Let $U$ be an infinite connected subset of $V$ with $\partial U \neq \emptyset$. The set $U$ is called a $D_p$-massive subset if there exists a nonnegative function $u \in BD_p(G)$ that satisfies the following:
\begin{enumerate} 
  \item []
   \begin{enumerate}  
     \item $\Delta_p u(x) = 0$ for all $x \in U,$
     \item $u(x) = 0$ for $x \in \partial U, $
     \item $\sup_{x \in U} u(x) = 1.$
   \end{enumerate} 
\end{enumerate} 
A function that satisfies the above conditions is called an {\em inner potential} of the $D_p$-massive subset $U$. The following will be needed in the proof of Lemma \ref{keyppara}.
\begin{Prop} \label{dpmassive}
If $U$ is a $D_p$-massive subset of $V$, then $\overline{i(U)}$ contains at least one point of $\partial_p(G)$.
\end{Prop}
\begin{proof}
We will write $\overline{U}$ for $\overline{i(U)}$, where the closure is taken in $Sp(BD_p(G))$. Assume $\overline{U} \cap \partial_p(G) = \emptyset$ and let $u$ be an inner potential for $U$. We may and do assume that $u = 0$ on $V \setminus U$. By the existence property for $p$-harmonic functions there exists a $p$-harmonic function $h_n$ on $B_n(o)$ such that $h_n = u$ on $\partial B_n(o)$ for each natural number $n$. Also $0 \leq \min_{y \in \partial B_n(o)} u(y) \leq h_n \leq \max_{y \in \partial B_n(o)} u(y) \leq 1$ on $B_n(o)$. Extend $h_n$ to all of $V$ by setting $h_n = u$ on $V \setminus B_n(o)$. By the minimizing property of $p$-harmonic functions, $I_p (h_n, B_n(o)) \leq I_p( u, B_n(o))$, consequently $I_p( h_n, V) \leq I_p( u, V)$. Both $h_n$ and $u$ are $p$-harmonic on $U \cap B_n(o)$ and $u(x) \leq h_n(x)$ for all $x \in \partial( U \cap B_n(o))$. The comparison principle says that $u \leq h_n$ on $U \cap B_n(o)$. On $B_n(o) \setminus U, u =0$, so $u \leq h_n \leq 1$ for each $n$. By taking a subsequence if needed we assume that $(h_n)$ converges pointwise to a function $h$. Now $u \leq h \leq 1$ on $V$, so $\sup_{x \in U} h(x) = 1$. By the convergence property for $p$-harmonic functions, $h$ is $p$-harmonic and $h \in BHD_p(G)$ since $I_p(h_n, V) \leq I_p(u, V) < \infty$ for all $n$.

Let $x \in \partial_p(G)$. Since $u - h_n = 0$ on $V \setminus B_n(o)$, we see that $\hat{u}(x) - \widehat{h_n} (x) = 0$ for all $n$, thus $\widehat{u - h} = 0$ on $\partial_p(G)$. According to Theorem \ref{vanishbound} $u - h \in B(\overline{\ell^p(G)})_{D_p}$. Hence $u = f + h$, where $f \in B(\overline{\ell^p(G)})_{D_p}$. Another appeal to Theorem \ref{vanishbound} shows that $\hat{u} = \hat{h}$ on $\partial_p(G)$. If $x \in \partial_p(G)$, then $\hat{u}(x) = 0$ because if $(x_n)$ is a sequence in $V$ converging to $x$, then $u(x_n) = 0$ for all but a finite number of $n$ since we are assuming $\overline{U} \cap \partial_p(G) = \emptyset$. So $\hat{h}(x) = 0$ for all $x \in \partial_p(G)$. Hence $h = 0$ on $V$ by the maximum principle, which contradicts $\sup_U h = 1$. Therefore, if $U$ is a $D_p$-massive subset of $V$, then $\overline{U}$ contains at least one point of $\partial_p(G).$
\end{proof}
It would be nice to know if the converse of the above proposition is true. That is, if $x \in \partial_p(G)$, then there exists a $D_p$-massive subset $U$ of $V$ such that $x \in \overline{U}$. The following result leads to a partial converse to Proposition \ref{dpmassive} and also describes a base of neighborhoods for open sets in $\partial_p(G)$. 
\begin{Prop} \label{base} 
Let $x \in \partial_p(G)$ and let $O$ be an open set in $\partial_p(G)$ containing $x$. Then there exists a subset $U$ of $V$ that satisfies the following:
\begin{enumerate} 
  \item []
    \begin{enumerate}  
      \item $U = \cup_{\alpha \in I}A_{\alpha}$, where each $A_{\alpha}$  is a $D_p$-massive subset of $V$ and $I$ is an index set. Moreover, if $\alpha \neq \beta$, then $A_{\alpha} \cap A_{\beta} = \emptyset.$
      \item $x \in \overline{U} \cap \partial_p(G) \subseteq O$.
    \end{enumerate} 
\end{enumerate}
\end{Prop}
\begin{proof}
Let $x \in \partial_p(G)$ and let $O$ be an open set of $\partial_p(G)$ containing $x$. By Urysohn's lemma there exists an $f \in C(Sp(BD_p(G)))$ with $0 \leq f \leq 1, f(x) = 1$ and $f = 0$ on $\partial_p(G) \setminus O$. Since the Gelfand transform of $BD_p(G)$ is dense in $C(Sp(BD_p(G)))$ with respect to the sup-norm we will assume $f \in BD_p(G)$. By Theorem \ref{decomp} we have the decomposition $f = w + h$, where $w \in B(\overline{\ell^p(G)})_{D_p}$ and $h \in BHD_p(G)$. Since $\hat{w} = 0$ on $\partial_p(G)$, it follows that $\hat{h}(x) = 1$ and $\hat{h} = 0$ on $\partial_p(G) \setminus O$. Moreover, $0 \leq \hat{h} \leq 1$ on $\partial_p(G)$ so $0 < h < 1$ on $V$ by the maximum principle and $0 \leq \hat{h} \leq 1$ on $Sp(BD_p(G))$ due to the density of $V$. Fix $\epsilon$ with $0 < \epsilon < 1$ and set $U = \{ x\in V \mid h(x) > \epsilon \}$. Let $A$ be a component of $U$. It now follows from the comparison principle that $A$ is infinite. Define a function $v$ on $V$ by 
$$ v = \frac{h- \epsilon}{1-\epsilon}. $$
There exists a $p$-harmonic function $u_n$ on $B_n(o) \cap A$ that takes the values $\mbox{max}\{0, v\}$ on $V\setminus (B_n(o) \cap A)$ such that $0 \leq u_n \leq 1$ on $B_n(o) \cap A$. By passing to a subsequence if necessary we may assume that the sequence $(u_n)$ converges pointwise to a function $u$. By the convergence property $u$ is $p$-harmonic on $A$. Also $v \leq u_n \leq 1$ on $B_n(o)$ so by replacing $u$ by a suitable scalar multiple if necessary we have $\sup_{a \in A} u(a) = 1$. Furthermore $u = 0$ on $\partial A$ because $h \leq \epsilon$ on $\partial A$. Since $h \in BD_p(G)$ it follows that $u \in BD_p(G)$. Thus $A$ is a $D_p$-massive subset with inner potential $u$. Hence, each component of $U$ is a $D_p$-massive subset in $V$. So $U = \cup_{\alpha \in I} A_{\alpha}$, where each $A_{\alpha}$  is $D_p$-massive. The proof of part (a) is complete.

Clearly $x \in \overline{U}$. We will now show that $\overline{U} \cap \partial_p(G) \subseteq O$. Let $y \in \overline{U} \cap \partial_p(G)$ and let $(y_k)$ be a sequence in $U$ that converges to $y$. Then $f(y) = \hat{h} (y) = \lim_{k \rightarrow \infty} h(y_k) \geq \epsilon$. Hence $y \in O$ since $f=0$ on $\partial_p(G) \setminus O.$
\end{proof}
The following partial converse to Proposition \ref{dpmassive} is a direct consequence of the preceding proposition.
\begin{Cor} \label{conversedpmassive}
If $\#(\partial_p(G))$ is finite, then for each $x \in \partial_p(G)$ there exists a $D_p$-massive subset $U$ of $V$ such that $x \in \overline{U}$.
\end{Cor}

\section{Proofs of Theorem \ref{pparaconstant} and Theorem \ref{extendcont}}\label{Proofsmain}
We begin this section by proving Theorem \ref{pparaconstant}. The key ingredient in the proof of the theorem is the following
\begin{Lem} \label{keyppara} 
Let $1 < p \in \mathbb{R}$ and suppose that $G$ is a $p$-parabolic graph. If $f$ is a nonconstant function in $BHD_p(G)$, then $\sup_V f > \limsup_{d(o,x) \rightarrow \infty} f$.
\end{Lem}
\begin{proof} 
Suppose that $\limsup_{d(o,x)\rightarrow \infty} f(x) = \sup_V f = M$. Since $f$ is nonconstant there exists an $\epsilon > 0$ such that the set $W= \{ x\in V \mid f(x) > M - \epsilon \}$ is a proper infinite subset of $V$. Let $U$ be a component of $W$. If $U$ is finite, then we can construct an unique $p$-harmonic function $w$ on $U$ which agrees with $f$ on $\partial U$. Since $f$ is $p$-harmonic we have $f = w$ on $U$ by uniqueness. But if $x \in U$, then $w(x) \leq \max_{y \in \partial U} f(y) \leq M - \epsilon < f(x)$, a contradiction. Thus $U$ is infinite. Now set $h = (f - M + \epsilon)/ \epsilon$. There exists a natural number $N$ such that $B_n(o) \cap U \neq \emptyset$ for $n > N$. For $n > N$ let $u_n$ be a $p$-harmonic function on $B_n(o) \cap U$ that takes the values $\max\{ 0, h\}$ on $V \setminus (B_n(o) \cap U)$. Note that $u_n \geq 0$. Since $h$ is $p$-harmonic on $B_n(o) \cap U$ it follows from the comparison principle that $h \leq u_n \leq 1$ on $B_n(o) \cap U$. By taking a subsequence if necessary we may assume that the sequence $(u_n)$ converges pointwise to a function $u$. By the convergence property $u$ is $p$-harmonic on $U$. If $x \in \partial U$, then $f(x) \leq M - \epsilon$. Consequently, $u_n(x) = 0$ for all $n$, which implies $u(x) = 0$. Thus $u = 0$ on $\partial U$. Since $\sup_U h = 1$ we see that $\sup_U u =1$. Using the minimizing property for $p$-harmonic functions it can be shown that $I_p( u_n, U \cap B_n(o) ) \leq I_p(\max\{0, h\}, U \cap B_n(o) )$ and it follows from this inequality that $I_p( u_n, U) \leq I_p( h, U)$. Hence $I_p (u, U) < \infty$ because $I_p( h, V) < \infty$. Thus $U$ is a $D_p$-massive subset of $V$.

By Proposition \ref{dpmassive} $\overline{U} \cap \partial_p(G) \neq \emptyset$, which contradicts Proposition \ref{noelement} since we are assuming $G$ is $p$-parabolic. Therefore, $\sup_V f > \limsup_{d(o,x) \rightarrow \infty} f$.
\end{proof}
We can now prove Theorem \ref{pparaconstant}. Let $h \in HBD_p(G)$ and suppose that $h$ is nonconstant. Since $h$ is bounded, $\sup_V h = B < \infty$. Lemma \ref{keyppara} says that there exists an $x \in V$ such that $h(x) = B$. By the maximum principle $h$ is constant on $V$, a contradiction. Hence $BHD_p(G)$ consists of only the constant functions. Therefore, $HD_p(G)$ is precisely the constant functions by Lemma 4.4 of \cite{H2}. This concludes the proof of Theorem \ref{pparaconstant}.

We now proceed to prove Theorem \ref{extendcont}. Suppose $f$ is a continuous function on $\partial_p(G)$. By Tietze's extension theorem there exists a continuous extension of $f$, which we also denote by $f$, to all of $Sp(BD_p(G))$. Let $(f_n)$ be a sequence in $BD_p(G)$ that converges to $f$ in the sup-norm. For each $n \in \mathbb{N}$ and each $r \in \mathbb{N}$ let $h_{n,r}$ be a function on $V$ that is $p$-harmonic on $B_r(o)$ and takes the values $f_n$ on $V\setminus B_r(o)$. The function $h_{n,r} \in BD_p(G)$ since $B_r(o)$ is finite and $\vert h_{n,r} \vert \leq \sup_V \vert f_n \vert$ because $\min_{y \in \partial B_r(o)} f_n(y) \leq h_{n,r} \leq \max_{y \in \partial B_r(o)} f_n(y)$ on $B_r(o)$. By the Ascoli-Arzela Theorem there exists a subsequence of $(h_{n,r})$, which we also denote by $(h_{n,r})$, that converges uniformly on all finite subsets of $V$ to a function $h_n$ as $r$ goes to infinity. The function $h_n$ is $p$-harmonic on $V$ by the convergence property. For each $r$ the minimizing property of $p$-harmonic functions gives $I_p(h_{n,r}, B_r(o)) \leq I_p(f_n, B_r(o))$, so $I_p(h_{n,r}, V) \leq I_p(f_n, V)$, which implies $h_n \in BHD_p(G)$.

Let $\epsilon > 0$. Since $(f_n) \rightarrow f$ in the sup-norm there exists a number $N$ such that for $n,m \geq N \sup_V \vert f_n - f_m \vert < \epsilon$. It follows that for all $r \in \mathbb{N}, \sup_{\partial B_r(o)} \vert h_{n,r} - h_{m,r} \vert < \epsilon$ because $f_n = h_{n,r}$ on $V\setminus B_r(o)$. Both $h_{n,r}$ and $h_{m,r} + \epsilon$ are $p$-harmonic on $B_r(o)$ and $h_{m,r} - \epsilon \leq h_{n,r} \leq h_{m,r} + \epsilon$ on $\partial B_r (o)$, so by applying the comparison principle we obtain $\sup_{B_r(o)} \vert h_{n,r} - h_{m,r} \vert < \epsilon$ for all $r$. It now follows that $\sup_{B_r(o)} \vert h_n - h_m \vert < 3 \epsilon$ for all $r$. Thus $\sup_V \vert h_n - h_m \vert \leq 3\epsilon$. Hence, the Cauchy sequence $(h_n)$ converges uniformly on finite subsets of $V$ to a function $h$, which is $p$-harmonic by the convergence property.

Let $\epsilon >0$, so there exists a number $N$ such that if $n \geq N$,
$$ \sup_V \vert f_n - f \vert < \epsilon \mbox{ and } \sup_V \vert h_n - h \vert < \epsilon.$$
Let $x \in \partial_p(G)$, since $f_n(x) = h_n(x)$ there exists a neighborhood $U$ of $x$ such that for all $y \in U, \vert h_n(y) - f_n(x) \vert < \epsilon$. Therefore, $\lim_{k \rightarrow \infty} h(x_k) = f(x)$, where $(x_k)$ is a sequence in $V$ that converges to $x$. Theorem \ref{extendcont} is now proven.

\section{Proofs of Theorem \ref{homeomorphic} and Theorem \ref{bijection}}\label{roughproofs}
Let $G$ and $H$ be graphs with vertex sets $V_G$ and $V_H$ respectively. Fix a vertex $o_G$ in $G$ and a vertex $o_H$ in $H$. Let $\phi \colon G \rightarrow H$ be a rough isometry and let $\phi^{\ast}$ denote the map from $\ell^{\infty}(H)$ to $\ell^{\infty}(G)$ given by $\phi^{\ast} f(x) = f (\phi(x))$. We start by defining a map $\bar{\phi} \colon \partial_p(G) \rightarrow \partial_p(H)$. Let $x \in \partial_p(G)$, then there exists a sequence $(x_n)$ in $V_G$ such that $(x_n) \rightarrow x$. Now $(\phi (x_n))$ is a sequence in the compact Hausdorff space $Sp (BD_p(H))$. By passing to a subsequence if necessary we may assume that $(\phi (x_n))$ converges to a unique limit $y$ in $Sp (BD_p(H))$. Now define $\bar{\phi} (x) = y$. Before we show that $y \in \partial_p (H)$ and $\bar{\phi}$ is well-defined we need the following lemma.
\begin{Lem} \label{adjoint}
Let $G$ and $H$ be graphs. If $\phi \colon G \rightarrow H$ is a rough isometry, then 
\begin{enumerate} 
  \item[]
  \begin{enumerate}
      \item $\phi^{\ast} \colon BD_p(H) \mapsto BD_p(G)$
      \item $\phi^{\ast} \colon \ell^p(H) \mapsto \ell^p(G) $
      \item $\phi^{\ast} \colon B(\overline{\ell^p(H)})_{D_p} \mapsto B(\overline{\ell^p(G)})_{D_p}.$
   \end{enumerate} 
 \end{enumerate}
\end{Lem}

\begin{proof}
We will only prove part (a) since the proofs of parts (b) and (c) are similar. Let $f \in BD_p(H)$. We will now show that $\phi^{\ast} f \in BD_p(G)$. Let $x \in V_G$ and $w \in N_x$, so $x$ and $w$ are neighbors in $G$ but $\phi(w)$ and $\phi(x)$ are not necessarily neighbors in $H$. However by the definition of rough isometry there exists constants $a \geq 1$ and $b \geq 0$ such that $d_H( \phi(w), \phi(x)) \leq a+b$. Set $h_1 = \phi (x)$ and $h_l = \phi(w)$ and let $ h_1, \dots, h_l $ be a path in $H$ with length at most $a + b$. Thus
\begin{equation}
\begin{split}
 |\phi^{\ast}f(w) - \phi^{\ast}f(x) |^p & =  | f(\phi(w)) - f(\phi(x)) |^p \\
  & \leq |a +b|^{p-1} \sum_{j=1}^{l-1} | f(h_{j+1}) - f(h_j)|^p. \label{eq:boundaplusb}
\end{split}
\end{equation}
The above inequality follows from Jensen's inequality applied to the function $x^p$ for $x >0$.

Let $y \in V_H$ and $z \in N_y$. We now claim that there is at most a finite number of paths in $H$ of length at most $a + b$ that contain the edge $y, z$ and have the endpoints $\phi(x)$ and $\phi(w)$. To see the claim let $U$ be the set of all elements in $V_G$ such that the following four distances: $d_H (\phi(x), y), d_H (\phi(x), z), d_H(\phi(w), y)$ and $d_H(\phi(w), z)$ are all at most $a + b$. Let $x$ and $x'$ be elements in $U$. By the triangle inequality, $d_H ( \phi(x'), \phi (x)) \leq d_H(\phi (x'), y) + d_H ( \phi(x), y)$. It now follows from the definition of rough isometry that $d_G (x', x) \leq 2a^2 + 3ab$. Thus the metric ball $B(x, 2a^2 + 3ab +1)$ contains $U$ as a subset. Hence the cardinality of $U$ is bounded above by some constant $k$. Observe that $k$ is independent of $y$ and $z$. Since $f \in BD_p(H)$ it follows from \ref{eq:boundaplusb} that
$$ \sum_{x \in V_G} \sum_{w \in N_x} | \phi^{\ast} f(w) - \phi^{\ast} f(x) |^p \\
\leq |a+b|^{p-1} k \sum_{y \in V_H} \sum_{z \in N_y} | f(z) - f(y) |^p \\
< \infty. $$
\end{proof}

We are now ready to prove 
\begin{Prop} \label{welldefined}
The map $\bar{\phi}$ is well-defined from $\partial_p(G)$ to $\partial_p(H)$.
\end{Prop}
\begin{proof}
Let $x, y$ and $(x_n)$ be as above. We start by showing that $y \in \partial_p(H)$. Lemma \ref{infinitelength} tells us that $d_G(o_G, x_n) \rightarrow \infty$ as $n \rightarrow \infty$. The element $\phi (o_G)$ is fixed in $H$ so it follows from the definition of rough isometry that $d_H ( \phi(o_G), \phi (x_n)) \rightarrow \infty$ as $n \rightarrow \infty$. Thus $ y \in Sp(BD_p(H)) \setminus H$ since $y = \lim_{n \rightarrow \infty} \phi (x_n) \notin H$. Let $f \in B(\overline{\ell^p(H)})_{D_p}$ and suppose $\hat{f}(y) \neq 0$. Then $0 \neq \lim_{n \rightarrow \infty} f( \phi (x_n)) = \phi^{\ast} f(x)$. By part (c) of Lemma \ref{adjoint} $\phi^{\ast} f \in B(\overline{\ell^p(G)})_{D_p}$ and Theorem \ref{vanishbound} says that $\phi^{\ast} f(x) = 0$, a contradiction. Hence $\hat{f}(y) =0 $ for all $f \in B(\overline{\ell^p(H)})_{D_p}$, so $y \in \partial_p (H)$.

We will now show that $\bar{\phi}$ is well-defined. Let $(x_n)$ and $(x_n')$ be sequences in $V_G$ that both converge to $x \in \partial_p(G)$. Now suppose that $(\phi (x_n))$ converges to $y_1$ and $(\phi (x_n'))$ converges to $y_2$ in $Sp( BD_p(H))$. Assume that $y_1 \neq y_2$ and let $f \in BD_p(H)$ such that $f(y_1) \neq f(y_2)$. By  part (a) of Lemma \ref{adjoint} $\phi^{\ast} f \in BD_p(G)$. Thus 
$$ \lim_{n \rightarrow \infty} \phi^{\ast} f(x_n) = \phi^{\ast} f(x) = \lim_{n \rightarrow \infty} \phi^{\ast} f(x_n')$$
which implies $f(y_1) = f(y_2)$, a contradiction. Hence $\bar{\phi}$ is a well-defined map from $\partial_p(G)$ to $\partial_p(H)$.
\end{proof}

The next lemma will be used to show that $\bar{\phi}$ is one-to-one and onto.
\begin{Lem} \label{ontolemma}
Let $\phi \colon G \rightarrow H$ be a rough isometry and let $\psi$ be a rough inverse for $\phi$. If $f \in D_p(G)$, then $\lim_{d_G(o_G,x) \rightarrow \infty} | f((\psi \circ \phi )(x)) - f(x) | = 0$.
\end{Lem}
\begin{proof}
Let $x \in V_G$, since $\psi$ is a rough inverse of $\phi$ there are non-negative constants $a, b$ and $c$ with $a \geq 1$ such that $d_G( (\psi \circ \phi) (x), x) \leq a(c + b)$. Let $x_1, x_2, \dots , x_n$ be a path in $V_G$ of length not more than $a (c+b)$ with $x_1 = x$ and $x_n = (\psi \circ \phi) (x)$. So
\begin{equation*}
\begin{split}
| f((\psi \circ \phi) (x)) - f(x) |^p & = |\sum_{k=1}^{n-1} ( f(x_{k+1}) - f(x_k))|^p \\
                              & \leq n^{p-1} \sum_{k = 1}^{n-1} | f(x_{k+1}) - f(x_k) |^p.
\end{split}
\end{equation*} 
The last sum approaches zero as $d_G(o_G, x) \rightarrow \infty$ since $f \in D_p(G)$ and $n \leq a(c+b)$. Thus $\lim_{d_G(o_G, x) \rightarrow \infty} | f( (\psi \circ \phi) (x)) - f(x) | = 0$.
\end{proof}

The next proposition shows that $\bar{\phi}$ is a bijection.
\begin{Prop}
The function $\bar{\phi}$ is a bijection.
\end{Prop}
\begin{proof}
Let $x_1, x_2 \in \partial_p(G)$ such that $x_1 \neq x_2$ and let $f \in BD_p(G)$ with $f(x_1) \neq f(x_2)$. There exists sequences $(x_n)$ and $(x_n')$ in $V_G$ such that $(x_n) \rightarrow x_1$ and $(x_n') \rightarrow x_2$. Now assume that $\bar{\phi}(x_1) = \lim_{n \rightarrow\infty} (\phi(x_n)) = \lim_{n \rightarrow \infty} ( \phi(x_n')) = \bar{\phi}(x_2)$, so $\lim_{n \rightarrow \infty} f( (\psi \circ \phi)(x_n)) = \lim_{n \rightarrow \infty} f ((\psi \circ \phi) (x_n'))$. It follows from Lemma \ref{ontolemma} that $\lim_{n \rightarrow \infty} f(x_n) = \lim_{n \rightarrow \infty} f(x_n')$, thus $f(x_1) = f(x_2)$, a contradiction. Hence $\bar{\phi}$ is one-to-one.

We now proceed to show that $\bar{\phi}$ is onto. Let $ y \in \partial_p(H)$ and let $(y_n)$ be a sequence in $V_H$ that converges to $y$. By passing to a subsequence if necessary, we can assume that there exist an unique $x$ in the compact Hausdorff space $Sp( BD_p(G))$ such that $( \psi(y_n)) \rightarrow x$. Since $\lim_{n \rightarrow \infty} d_H(o_H, y_n) \rightarrow \infty$ it follows that $\lim_{n \rightarrow \infty} d_G(o_G,\psi (y_n)) \rightarrow \infty$, so $x \notin G$. Using an argument similar to the first paragraph in the proof of Proposition \ref{welldefined} we obtain $x \in \partial_p(G)$. The proof will be complete once we show that $\bar{\phi} (x) = y$. Let $f \in BD_p(H)$. By Lemma \ref{ontolemma} we see that $\lim_{n \rightarrow \infty} |f( (\phi \circ \psi) (y_n)) - f(y_n) | = 0$. Thus $f(\bar{\phi} (x) ) = f(y)$ for all $f \in BD_p(H)$. Hence $\bar{\phi} (x) = y$.
\end{proof}

To finish the proof that the bijection $\bar{\phi}$ is a homeomorphism we only need to show that $\bar{\phi}$ is continuous, since both $Sp(BD_p(G))$ and $Sp(BD_p (H))$ are compact Hausdorff spaces. Let $W$ be an open set in $\partial_p (H)$ and let $x \in \bar{\phi}^{-1}(W)$. Let $y \in W$ for which $x = \bar{\phi}^{-1}(y)$. By Proposition \ref{base} there exists a subset $U$ of $V_H$ such that $y \in \overline{U}$ and $\overline{U} \cap \partial_p(H) \subseteq W$. Moreover, we saw in the proof of Proposition \ref{base} that there is an $h \in BHD_p(H)$ for which $\hat{h}(y) = 1, \hat{h} = 0$ on $\partial_p(H) \setminus W$ and $\hat{h} \geq \epsilon$ on $\overline{U}$, where $0 < \epsilon < 1$. By Lemma \ref{adjoint} (a) $\phi^{\ast} h = h \circ \phi \in BD_p(G)$. Combining Theorem \ref{decomp} and Theorem \ref{vanishbound} we have an $\bar{h} \in BHD_p(G)$ that satisfies $\bar{h} = \hat{h} \circ \bar{\phi}$ on $\partial_p(G)$. Let $O = \{ x' \in \partial_p(G) \mid \bar{h}(x') > \epsilon \}$. Now $O$ is an open set containing $x$ since $\bar{h}$ is continuous on $\partial_p(G)$ and $\bar{h}(x) = 1$. For $z \in O$ we see that $\hat{h}(\bar{\phi} (z)) = \bar{h}(z) \geq \epsilon$, thus $\bar{\phi}(z) \in W$ for all $z$ in $O$. Thus $O \subseteq \bar{\phi}^{-1} (W)$. Since our choice of $x$ was arbitrary, $\bar{\phi}^{-1} (W)$ is open and consequently $\bar{\phi}$ is continuous. The proof that $\bar{\phi}$ is a homeomorphism is complete.

We will now prove Theorem \ref{bijection}. Let $\phi$ be a rough isometry from $G$ to $H$ and let $\psi$ be a rough inverse for $\phi$. Let $h \in BHD_p(G)$. By part (a) of Lemma \ref{adjoint}, $h \circ \psi \in BD_p(H)$. Let $\pi(h \circ \psi)$ be the unique element in $BHD_p(H)$ given by Proposition \ref{decomp}. We now define a map $\Phi \colon BHD_p (G) \mapsto BHD_p(H)$ by $\Phi (h) = \pi(h \circ \psi).$ Theorem \ref{vanishbound} implies that $\pi (h \circ \psi)(\bar{\phi} (x)) = (h \circ \psi)(\bar{\phi}(x))$ for all $x \in \partial_p(G)$, where $\bar{\phi}$ is the homeomorphism from $\partial_p(G)$ to $\partial_p(H)$ defined earlier in this section. Thus $\Phi (h) (\bar{\phi}(x)) = (h \circ \psi)(\bar{\phi}(x)) = h(x)$ for all $x \in \partial_p(G)$. We can now show that $\Phi$ is one-to-one. Let $h_1, h_2 \in BHD_p(G)$ and suppose that $\Phi (h_1) = \Phi (h_2)$. So $\Phi (h_1) (\bar{\phi}(x)) = \Phi (h_2) ( \bar{\phi}(x))$ for all $ x \in \partial_p (G)$, which implies $h_1(x) = h_2(x)$ for all $x \in \partial_p (G)$. Hence, $h_1 = h_2$ by Proposition \ref{boundval}. Thus $\Phi$ is one-to-one.

We will now show that $\Phi$ is onto. Let $f \in BHD_p(H)$. Then $f \circ \phi \in BD_p(G)$. Let $h = \pi(f \circ \phi)$, where $\pi (f \circ \phi)$ is the unique element in $BHD_p(G)$ given by Proposition \ref{decomp}. Let $y \in \partial_p (H)$. Since $h(x) = \pi( f \circ \phi)(x)$ for all $ x \in \partial_p(G)$ and $\bar{\psi} \circ \bar{\phi}$ equals the identity on $\partial_p(G)$, we see that $(\Phi (h))(y) = \pi (h \circ \psi)(y) = h (\psi(y)) = f((\phi \circ \psi)(y)) = f(y)$. Thus $\Phi$ is onto and the proof of Theorem \ref{bijection} is complete.

The map $\Phi$ is an isomorphism in the case $p=2$ since $BHD_2 (G)$ and $BHD_2 (H)$ are linear spaces. However, in general these spaces are not linear if $ p \neq 2$.

\section{The first reduced $\ell^p$-cohomology of $\Gamma$}\label{lpcohomology}
In the last two sections of this paper $\Gamma$ will denote a finitely generated group with generating set $S$. So for a real-valued function $f$ on $\Gamma$ the $p$-th power of the gradient and the $p$-Laplacian of $x \in \Gamma$ are
\begin{equation*}
 \begin{split}
  \vert Df (x) \vert^p & =  \sum_{s \in S} \vert f(xs^{-1}) - f(x) \vert^p, \\
  \Delta_p f(x)          & =  \sum_{s \in S} \vert f(xs^{-1}) - f(x) \vert^{p-2} ( f(xs^{-1}) - f(x)).
 \end{split}
\end{equation*} 
If $f \in D_p(\Gamma)$, then
 $$ \left( \parallel f \parallel_{D_p} = I_p (f, \Gamma) + \vert f(e) \vert^p \right)^{\frac{1}{p}},$$
where $e$ is the identity element of $\Gamma$. Also $\ell^p (\Gamma)$ is the set that consists of real-valued functions on $\Gamma$ for which $\sum_{x \in \Gamma} \vert f(x) \vert^p$ is finite. The first reduced $\ell^p$-cohomology space of $\Gamma$ is defined by 
$$ \bar{H}^1_{(p)} (\Gamma) = D_p(\Gamma) / (\overline{\ell^p(\Gamma) \oplus \mathbb{R}})_{D_p}.$$
We now prove Theorem \ref{lpchar}. Suppose $\partial_p (\Gamma) = \emptyset$. By Proposition \ref{noelement} there exists a sequence $(f_n)$ in $\mathbb{R}\Gamma$ that satisfies $\parallel f_n - 1_{\Gamma} \parallel_{D_p} \rightarrow 0$. It follows that $I_p( f_n, \Gamma) \rightarrow 0$ and $(f_n(e)) \not\rightarrow 0$. Thus $\bar{H}^1_{(p)} (\Gamma) = 0$ by Theorem 3.2 of \cite{PulsCMB}. We now assume $\partial_p(G) \neq \emptyset$. It was shown in \cite[Theorem 3.5]{Puls} that $\bar{H}_{(p)}^1 (\Gamma) \neq 0$ if and only if $HD_p(\Gamma) \neq \mathbb{R}$. Since $\#(S) < \infty$, Lemma 4.4 of \cite{H2} says that $BHD_p(\Gamma) = \mathbb{R}$ if and only if $HD_p(\Gamma) =\mathbb{R}$. Theorem \ref{lpchar} now follows from Theorem \ref{plious}.

We now use Theorem \ref{lpchar} to compute $\partial_p(\Gamma)$ and $R_p(\Gamma)$ for some special cases of $\Gamma$. By \cite[Corollary 1.10]{HolaSoarliouville}, $BHD_p(\Gamma) = \mathbb{R}$ when $\Gamma$ has polynomial growth and $1 < p \in \mathbb{R}$. Thus, if $\Gamma$ has polynomial growth then $\bar{H}_{(p)}^1 (\Gamma) = 0$ and $\partial_p(\Gamma)$ is either the empty set or contains exactly one element. It would be nice to know when a group with polynomial growth is $p$-parabolic or $p$-hyperbolic. This has been worked out for the case $\Gamma = \mathbb{Z}^n$, where $n$ is a positive integer. Example 4.1 in \cite{Yamasaki1977} showed that $\mathbb{Z}$ is $p$-parabolic for $p > 1$, thus $\partial_p(\mathbb{Z}) = \emptyset$ for $p > 1$. The main result of \cite{Maeda} says that for $n \geq 2, \mathbb{Z}^n$ is $p$-parabolic if and only if $p \geq n$. Hence, $\partial_p(\mathbb{Z}^n) = \emptyset$ if $p \geq n$ and $\partial_p(\mathbb{Z}^n)$ consists of exactly one point if $1 < p < n$.

There is a one-to-one correspondence between the maximal ideals of $BD_p(\Gamma)$ and the points of $Sp (BD_p(\Gamma))$. If $\tau \in R_p(\Gamma)$, then $ker(\tau)$ is the maximal ideal of $BD_p(\Gamma)$ that corresponds to $\tau$. For each $x \in \Gamma, \delta_x \in ker(\tau)$. Due to the continuity of $\tau$ we see that $\ell^p(\Gamma) \subseteq ker(\tau)$. Now assume that $\Gamma$ is nonamenable. Then by \cite[Corollary 1]{guichardet} $\ell^p(\Gamma)$ is closed in $D_p(\Gamma)$. Hence $(\overline{\mathbb{R}\Gamma})_{D_p} = \ell^p(\Gamma)$. Furthermore, $(\overline{\ell^p(\Gamma)})_{BD_p} = \ell^p(\Gamma)$ because $(\overline{\ell^p(\Gamma)})_{BD_p} \subseteq B(\overline{\ell^p(\Gamma)})_{D_p}$. Thus $\hat{f}(\tau) = 0$ for every $f \in (\overline{\mathbb{R}\Gamma})_{D_p}$. Therefore, $R_p(\Gamma) = \partial_p(\Gamma)$ when $\Gamma$ is nonamenable. Consequently, $R_p(\Gamma)$ contains exactly one point when $\Gamma$ is nonamenable and $\bar{H}_{(p)}^1 (\Gamma) = 0$. Examples of groups that satisfy this last condition for $1 < p < \mathbb{R}$ are nonamenable groups with infinite center \cite[Theorem 4.2]{MartinValette}; $\Gamma_1 \times \Gamma_2 \times \cdots \times \Gamma_n$ where $n \geq 2$, each $\Gamma_i$ is finitely generated and at least one of the $\Gamma_i$ is nonamenable \cite[Theorem 4.7]{MartinValette}.

\section{Translation Invariant Linear Functionals}\label{tilf}
Recall that $\Gamma$ denotes a finitely generated group with generating set $S$. In this section we will study $TILFs$ on $D_p(\Gamma)/\mathbb{R}$. By definition we have the following inclusions:
$$ \mbox {Diff}(\ell^p(\Gamma)) \subseteq \mbox{Diff}(D_p(\Gamma)/\mathbb{R}) \subseteq \ell^p(\Gamma) \subseteq D_p(\Gamma)/\mathbb{R}.$$
The set $D_p(\Gamma)/\mathbb{R}$ is a Banach space under the norm induced from $I_p(\cdot, \Gamma)$. Thus if $[f]$ if a class from $D_p(G) / \mathbb{R}$, then its norm is given by 
$$ \parallel [f] \parallel_{D(p)} = \left( \sum_{x \in \Gamma} \sum_{s \in S} \vert f(xs^{-1}) - f(x) \vert^p \right)^{1/p}.$$
We will write $\parallel f \parallel_{D(p)}$ for $\parallel [f] \parallel_{D(p)}$. Now $(\overline{\ell^p (\Gamma)})_{D(p)} = D_p(\Gamma)/\mathbb{R}$ if and only if $( \overline{ \ell^p(\Gamma) \oplus \mathbb{R}})_{D_p} = D_p(\Gamma).$ So $\bar{H}_{(p)}^1 (\Gamma) = 0$ if and only if $( \overline{\ell^p (\Gamma)})_{D(p)} = D_p(\Gamma)/\mathbb{R}$. We begin by proving the following
\begin{Lem} \label{equality}
$(\overline{\mbox{Diff}(D_p(\Gamma)/\mathbb{R})})_{D(p)} = (\overline{\ell^p(\Gamma)})_{D(p)}$.
\end{Lem}
\begin{proof}
Let $f \in \ell^p(\Gamma)$. By \cite[Lemma 1]{Woodward} there is a sequence $( f_n )$ in $\mbox{Diff}(\ell^p(\Gamma))$ that converges to $f$ in the $\ell^p$-norm. It follows from Minkowski's inequality that for $s \in S,\parallel (f-f_n)_s - (f-f_n) \parallel_p^p = \sum_{x \in \Gamma} | f(xs^{-1}) - f_n (xs^{-1}) - (f(x) - f_n(x))|^p \rightarrow 0$ as $n \rightarrow \infty$. Hence $f \in (\overline{\mbox{Diff}(\ell^p(\Gamma))})_{D(p)}$ which implies $\ell^p(\Gamma) \subseteq ( \overline{\mbox{Diff}(\ell^p(\Gamma))})_{D(p)}$. The result now follows.
\end{proof}

We are now ready to prove the following characterization, which is a direct consequence of the Hahn-Banach theorem, of nonzero continuous TILFs on $D_p(\Gamma)/\mathbb{R}$.
\begin{Thm} \label{charlptilf}
Let $1 < p \in \mathbb{R}$. Then $\bar{H}^1_{(p)}(\Gamma) \neq 0$ if and only if there exists a nonzero continuous TILF on $D_p(\Gamma)/\mathbb{R}$.
\end{Thm}
\begin{proof}
If $\bar{H}^1_{(p)}(\Gamma) \neq 0$, then $(\overline{\ell^p(\Gamma)})_{D(p)} \neq D_p(\Gamma)/\mathbb{R}$. It now follows from the Hahn-Banach theorem that there exists a nonzero continuous linear functional $T$ on $D_p(\Gamma)/\mathbb{R}$ such that $(\overline{\ell^p(\Gamma)})_{D(p)}$ is contained in the kernel of $T$. Thus $T$ is translation invariant by Lemma \ref{equality}.

Conversely if $T$ is a continuous TILF on $D_p(\Gamma)/\mathbb{R}$, then $T(f) = 0$ for all $f \in (\overline{\ell^p(\Gamma)})_{D(p)}$. So if there exists a nonzero continuous TILF on $D_p(\Gamma)/\mathbb{R}$, then $(\overline{\ell^p(\Gamma)})_{D(p)} \neq D_p(\Gamma)/\mathbb{R}$.
\end{proof}

Theorem \ref{chartilf} now follows by combining Theorems \ref{charlptilf} and \ref{lpchar}.

If $h \in D_p(\Gamma)/\mathbb{R}$, then $\langle \triangle_p h, \cdot \rangle$ is a well-defined continuous linear functional on $D_p(\Gamma)/\mathbb{R}$ since equivalent functions in $D_p(\Gamma)/\mathbb{R}$ differ by a constant. It was shown in Proposition 3.4 of \cite{Puls} that if $h \in HD_p(\Gamma)/ \mathbb{R}$ and $f \in ( \overline{\ell^p(\Gamma)})_{D(p)}$, then $\langle \triangle_p h, f \rangle = 0$. Consequently, if $h \in HD_p(\Gamma)/\mathbb{R}$, then $\langle \triangle_p h, \cdot \rangle$ defines a continuous TILF on $D_p(\Gamma)/ \mathbb{R}$. Thus there are no nonzero continuous TILFs on $D_p(\Gamma)/ \mathbb{R}$ when $HD_p(\Gamma)$ only contains the constant functions.

If $\bar{H}_{(p)}^1(\Gamma) = 0$, then $(\overline{\ell^p(\Gamma)})_{D(p)} = D_p(\Gamma)/\mathbb{R}$. It is known that $\ell^p(\Gamma)$ is closed in $D_p(\Gamma)/\mathbb{R}$ if and only if $\Gamma$ is nonamenable, \cite[Corollary 1]{guichardet}. As was mentioned in Section \ref{Outlinestatementmain}, if $\Gamma$ is nonamenable, then zero is the only TILF on $\ell^p(\Gamma)$. Consequently zero is the only TILF on $D_p(\Gamma)/\mathbb{R}$ when $\Gamma$ is nonamenable and $\bar{H}_{(p)}^1 (\Gamma) =0$. Summing up we obtain:
\begin{Thm}
Let $\Gamma$ be an infinite, finitely generated group and let $1 < p \in \mathbb{R}$. The following are equivalent
\begin{enumerate}
\item $\bar{H}^1_{(p)} (\Gamma) = 0$
\item $\mbox{Either } \partial_p(\Gamma) = \emptyset \mbox{ or } \#(\partial_p(\Gamma)) = 1$
\item $HD_p(\Gamma) = \mathbb{R}$
\item $BHD_p(\Gamma) = \mathbb{R}$
\item The only continuous TILF on $D_p(\Gamma)/\mathbb{R}$ is zero\\
If moreover $\Gamma$ is nonamenable, then this is still equivalent to:
\item Zero is the only TILF on $D_p(\Gamma)/\mathbb{R}$
\end{enumerate}
\end{Thm}

We will now give some examples that show zero is not the only $TILF$ on $D_p(\Gamma)/\mathbb{R}$ when $\Gamma$ is nonamenable, this differs from the $\ell^p(\Gamma)$ case. It was shown in \cite[Corollary 4.3]{Puls} that $\bar{H}^1_{(p)} (\Gamma) \neq 0$ for groups with infinitely many ends and $1 < p \in \mathbb{R}$. Thus by Theorem \ref{charlptilf} there exists a nonzero continuous TILF on $D_p(\Gamma)/\mathbb{R}$. Another question that now arises is: if there is a nonzero continuous TILF on $D_r(\Gamma)/\mathbb{R}$ for some nonamenable group $\Gamma$ and some real number $r$, then is it true that there is a nonzero continuous TILF on $D_p(\Gamma)/\mathbb{R}$ for all real numbers $p>1$? The answer to this question is no. To see this let $\mathcal{H}^n$ denote hyperbolic $n$-space, and suppose $\Gamma$ is a group that acts properly discontinuously on $\mathcal{H}^n$ by isometries and that the action is cocompact and free. By combining Theorem 2 of \cite{BMV} and Theorem 1.1 of \cite{PulsADM} we obtain $\bar{H}^1_{(p)}(\Gamma) \neq 0$ if and only if $p > n-1$.

\bibliographystyle{plain}
\bibliography{revisegraphpharmbound}
\end{document}